\documentclass[12pt]{article}
\usepackage[mathscr]{eucal}
\usepackage{amssymb,latexsym}
\usepackage{verbatim}
\usepackage{amsmath}
\usepackage{amsthm}
\usepackage{enumerate}
\usepackage{authblk}
\usepackage{color}
\usepackage{hyperref}

\newtheorem{thm}{Theorem}[section]
\newtheorem{prop}[thm]{Proposition}%[section]
\theoremstyle{remark} %added by Hyun
\newtheorem{rem}{Remark}%[section]
\theoremstyle{plain} %added by Hyun
\renewcommand\l{\lambda}

\newcommand\bbR{{\mathbb R}}

\newcommand\vep{\varepsilon}

\renewcommand\S{\Sigma}

\renewcommand\d{\partial}

\newcommand\f{\phi}

\newcommand\D{\nabla}

\renewcommand\div{{\rm div}}

\newcommand\<{\langle}
\renewcommand\>{\rangle}

\renewcommand\l{\lambda}

\renewcommand\th{\theta}

\newcommand\beq{\begin{equation}}
\newcommand\eeq{\end{equation}}
\newcommand\ben{\begin{enumerate}}
\newcommand\een{\end{enumerate}}
\newcommand\bit{\begin{itemize}}
\newcommand\eit{\end{itemize}}

\DeclareMathOperator{\diver}{div}

\renewcommand{\div}{\diver}

\newcounter{mnotecount}

\title{Some scalar curvature warped product \\splitting theorems}    

\author[*]{Gregory J. Galloway}
\author[$\dag$]{Hyun Chul Jang}

\affil[*]{\small Department of Mathematics, 

University of Miami, Coral Gables, FL }

\affil[$\dag$]{Department of Mathematics,

University of Connecticut, Storrs, CT}

\begin{document}
\date{}
\maketitle
\vspace{.2in}

\begin{abstract}  We present several rigidity results for Riemannian manifolds $(M,g)$ with scalar curvature $S \ge -n(n-1)$ (or $S\ge 0$), and having compact boundary $N$ satisfying a related mean curvature inequality.  The proofs make use of results on marginally outer trapped surfaces applied to appropriate initial data sets.  One of the results involves an analysis of Obata's equation on manifolds with boundary. This result is relevant to recent work of  Lan-Hsuan Huang and the second author concerning the rigidity of asymptotically locally hyperbolic manifolds with zero mass.

\end{abstract}

\section{Introduction}
\label{intro}

In the recent development of differential geometry, it has been of increasing interest to understand the geometry of Riemannian manifolds with lower bounds on their scalar curvature. In \cite{sysc,sysc2}, R. Schoen and S.-T. Yau proved the milestone result that the $n$-dimensional torus $T^n$, for $3\le n\le 7$, does not admit a metric of positive scalar curvature  by using minimal surface techniques.   In more recent work,  Schoen and Yau \cite{SY3} have been able to use the minimal surface method to prove this for all dimensions $n \ge 3$.   This had been proved by M. Gromov and H. B. Lawson \cite{GL} using spinor methods.
The key observation made by Schoen and Yau in \cite{sysc2} 
is the following.
\begin{prop}\label{SYstable}
	Let $(M^n,g),n\ge 3$ be a Riemannian manifold with positive scalar curvature, $S  >0$. If $N^{n-1}$ is a stable, two-sided closed minimal hypersurface in $M^n$ then $N^{n-1}$ admits a metric of positive scalar curvature.
\end{prop}
Moreover, by refinements of the arguments in \cite{sysc2}, one obtains the rigidity statement that  if $S \ge 0$, and $N^{n-1}$ does not admit a metric of positive scalar curvature then $N^{n-1}$ is totally geodesic and Ricci flat, and $S=0$ along $N^{n-1}$ (cf. \cite{FCS, Gschoen}).  In \cite{Cai}, M. Cai proved the following splitting theorem by assuming $N$ is area-minimizing instead of being only stable (see also \cite{Gschoen} for a simplified proof).
\begin{prop}\label{Caisplit}
	Let $(M^n,g),n\ge 3$ be a Riemannian manifold with nonnegative scalar curvature, $S\ge 0$, and suppose $N^{n-1}$ is a two-sided closed minimal hypersurface which locally minimizes area. If $N$ does not admit a metric of positive scalar curvature then there exists a neighborhood $V$ of $N$ such that $(V,g|_V)$ is isometric to $(-\delta,\delta)\times N$ with product metric $dt^2+h$, where $h=g|_N$, and $(N,h)$ is Ricci flat.
\end{prop}

This result extends to higher dimensions the torus splitting result in \cite{CG} for $3$-manifolds of nonnegative scalar curvature.   For some related rigidity results in three dimensions under different assumptions on the ambient scalar curvature and the topology of the minimal surface, see for example, \cite{Bray, Nunes, MicallefMoraru, Ambrozio, Eichmair1, Eichmair2}.

The minimal surface techniques introduced in \cite{sysc,sysc2} also played an important role for the proof of the celebrated positive mass theorem for asymptotically flat manifolds by Schoen and Yau in \cite{SY1,SY2}, which they have now extended to arbitrary dimension $n \ge 3$ in \cite{SY3}.   These results include the rigidity statement that the mass vanishes if and only if the manifold is isometric to  Euclidean space. Somewhat more relevant for the present work are results concerning asymptotically hyperbolic manifolds.  A proof of the positivity of mass in this setting was obtained by X. Wang \cite{Wang} for spin manifolds, with improvements by P. Chru\'sciel and M. Herzlich \cite{ChH}.   In the paper \cite{ACG}, L. Andersson, M. Cai, and the first author proved a positive mass result without spin assumption in dimensions $n$, $3 \le n \le 7$, for asymptotically hyperbolic manifolds, assuming a sign on the mass aspect.   As an element in the proof, a splitting result analogous to Proposition  \ref{Caisplit} was obtained in \cite[Section 2.2]{ACG}, whereby the `brane' functional takes the place of the area functional and the scalar curvature satisfies $S\ge -n(n-1)$.  Recently, making use of work of Lohkamp \cite{Loh}, Chru\'sciel and Delay in \cite{Delay:2019uy} have established 
the nonnegativity of the mass for asymptotically hyperbolic manifolds, without spin assumption and in arbitrary dimension $n \ge 3$.  The rigidity statement, when the mass vanishes, has been proved by L.-H. Huang, D. Martin and the second author in \cite{Huang:2019tm}. 

The aim of the present paper is to obtain splitting theorems for manifolds with compact boundary satisfying the scalar curvature inequality, $S\ge -\vep n(n-1)$ (where $\vep=0$ or $1$). An initial motivation for this paper comes from recent work of L.-H. Huang and the second author \cite{LanHyun} concerning the rigidity of asymptotically locally hyperbolic manifolds of zero mass.  

For our splitting results, we will use a condition that replaces the least area (or brane minimization) assumption.  Let $(M,g)$ be a Riemannian manifold with compact boundary $N$ having mean curvature $H_N \le H_0$, $H_0 \in \bbR$.\footnote{For simplicity we always assume $M$ and  $N$ are connected.} To set sign conventions, the mean curvature $H_N$ is defined as the divergence of the inward pointing unit normal.
We say that $N$ is {\it weakly outermost} if there does not exist a compact hypersurface $\S \subset M \setminus N$ cobordant to $N$ satisfying the (strict) mean curvature inequality, $H_{\S} < H_0$. We further define, in order to state the local version of our results, that $N$ is {\it locally weakly outermost} provided that there is a neighborhood $U$ of $N$ such that $N$ is weakly outermost in $(U,g|_U)$. 
In our results, in addition to a weakly outermost condition, we will also require that the boundary $N$ not admit a metric of positive scalar curvature, as in Proposition \ref{Caisplit}. We will discuss the necessity of these assumptions in Remark \ref{examples}.

We now state our main result.

\begin{thm}\label{splitlocal1}
Let $(M,g)$ be an $n$-dimensional ($n \ge 3$) Riemannian manifold with compact boundary $N$.  Assume:
\ben
\item $M$ has scalar curvature $S \ge  -\varepsilon n(n-1)$, where $\varepsilon = 0$ or $1$.
\item $N$ has mean curvature $H_N \le  \varepsilon (n-1)$.
\item $N$ does not carry a metric of positive scalar curvature and is locally weakly outermost. 
\een
Then there exists a neighborhood $V$ of $N$ such that 
$(V,g|_V)$ is isometric to $[0,\delta) \times N$, with (warped) product metric 
$dt^2 + e^{2\varepsilon t} h$,  where $(N,h)$ is Ricci flat. 
\end{thm}
If we assume $N$ is (globally) weakly outermost, one can obtain the global splitting result, as stated in Theorem \ref{split}. In the case  $\vep = 0$, the conclusion is that a neighborhood $V$ of $N$ splits as a product, which can be viewed as a variation of Proposition \ref{Caisplit}. On the other hand, in the case $\vep = 1$, $V$ splits as a warped product. Note that if $h$ is flat, the manifold $([0,\infty)\times N,dt^2+e^{2t} h)$ is of constant sectional curvature $-1$, and serves as a model space to define an asymptotically locally hyperbolic manifold.

\begin{rem}\label{examples}
	The assumption in point 3 that the boundary $N$ is weakly outermost is not sufficient to obtain the desired rigidity. For the case $\vep=0$, consider the {\it spatial Schwarzschild manifold}: $M = \bbR^n \setminus \{r < (\frac{m}{2})^\frac1{n-2}\}$, with metric (in isotropic coordinates), 
	$$
		g = \left(1+ \frac{m}{2r^{n-2}}\right)^\frac{4}{n-2} g_E \,,
	$$
	where $g_E$ is the Euclidean metric and $r = \sqrt{\sum_{i= 1}^n x_i^2}$.  $M$ has vanishing scalar curvature, $S =0$, and the boundary $N: r = 
	(\frac{m}{2})^\frac1{n-2}$ is minimal, $H_N = 0$.  Moreover, it follows from the maximum principle for hypersurfaces that $N$ is weakly outermost. However, the conclusion of Theorem \ref{splitlocal1} does not hold

	For the case $\vep=1$, the {\it AdS Schwarzschild manifold} is a further example illustrating the need for the scalar curvature assumption on $N$: $M=[r_m,\infty)\times S^{n-1}$ with metric
	\[
		g=\left(1+r^2-\frac{2m}{r^{n-2}}\right)^{-1}dr^2+r^2 g_{S^{n-1}}
	\]
	where $r_m=(2m)^{\frac{1}{n-2}}$ and $g_{S^{n-1}}$ is the standard unit sphere metric. In this case, $(M,g)$ has constant scalar curvature $S=-n(n-1)$ and the mean curvature of its boundary $N=\{r_m\}\times S^{n-1}$ is equal to $n-1$. Also, $N$ is weakly outermost but $N$ carries a metric of positive scalar curvature.

As a final example, consider the toroidal Kottler metrics with $m>0$: 
$M=[r_0,\infty)\times T^{n-1}$ with metric
	\[
		g=\left(r^2-\frac{2m}{r^{n-2}}\right)^{-1}dr^2+r^2 h
	\]
	where $r_0=(2m)^{\frac{1}{n}}$ and $h$ is a flat metric on $T^{n-1}$. One can easily check this example satisfies the conditions 1, 2 (with $\vep=1$).  $T^{n-1}$ does not carry a metric of positive scalar curvature, but the boundary $N=\{r_0\}\times T^{n-1}$ is not weakly outermost. This example shows that the boundary $N$ being weakly outermost is needed as well, in addition to $N$ not admitting a metric of positive scalar curvature.
\end{rem}

By similar arguments one also obtains the following variation of Theorem \ref{splitlocal1}.

\begin{thm}\label{splitlocal2}
Let $(M,g)$ be an $n$-dimensional ($n \ge 3$) Riemannian manifold with compact boundary $N$.  Assume:
\ben
\item $M$ has scalar curvature $S \ge  -\varepsilon n(n-1)$, where $\varepsilon = 0$ or $1$.
\item $N$ has mean curvature $H_N \le  -\varepsilon (n-1)$.
\item $N$ does not carry a metric of positive scalar curvature and is locally weakly outermost. 
\een
Then there exists a neighborhood $V$ of $N$ such that 
$(V,g|_V)$ is isometric to $[0,\delta) \times N$, with (warped) product metric 
$dt^2 + e^{-2\varepsilon t} h$,  where $(N,h)$ is Ricci flat. 
\end{thm}

The above variation can  be roughly regarded as a scalar curvature version of a (warped product) splitting theorem proved by Croke and Kleiner in \cite{Croke}, in which they assume the corresponding lower bound on Ricci curvature, but do not require a scalar curvature condition on the boundary $N$.  In this Ricci curvature case, the condition of being weakly outermost is implicit in their assumptions.

In addition to the above mentioned results, we prove a global splitting result by using Obata's equation, $\nabla^2 f=fg$. 

\begin{thm}\label{splithess}
	Let $(M,g)$ be an $n$-dimensional $(n\ge 3)$ complete, noncompact Riemannian manifold with compact boundary $ N$. Let $h=g|_ N$. Suppose that
	\begin{enumerate}
		\item $S\ge -n(n-1)$ in a neighborhood of $N$.
		\item $ N$ has mean curvature $H_ N\le \delta (n-1)$, where $\delta=1$ or $-1$.
		\item $ N$ does not carry a metric of positive scalar curvature and is locally weakly outermost.
		\item There exists a nonzero function $f$ satisfying $\nabla^2 f=fg$.
	\end{enumerate}
	Then $(M,g)$ is isometric to $[0,\infty)\times N$, with warped product metric $dt^2+e^{2\delta t} h$ where $( N,h)$ is Ricci flat.
\end{thm}
Here we only require $N$ to be locally weakly outermost. Instead, we can extend the local splitting result globally by assuming the existence of a nontrivial solution to Obata's equation. Note that the resulting warped product corresponds to an unbounded portion of the hyperbolic cusp: it contains either an expanding end when $\delta=1$ or a shrinking end when $\delta=-1$.
This result plays a role in the recent work of Lan-Hsuan Huang and the second author \cite{LanHyun} mentioned above.

As discussed in the next section, which includes relevant background, the proofs of 
Theorems~\ref{splitlocal1} and \ref{splitlocal2} make use of results on {\it marginally outer trapped surfaces}, applied to specific {\it initial data  sets}.  The proof of Theorem~\ref{splitlocal1}, and its globalization are presented in Section~3.  The proof of Theorem \ref{splithess} is presented in Section 4.

\section{Marginally outer trapped surfaces}

For the proof of Theorem \ref{splitlocal1}, we will make use of the theory of marginally outer trapped surfaces. Such surfaces play an important role in the theory of black holes, and, as indicated below, may be viewed as spacetime analogues of minimal surfaces in Riemannian geometry.  For further background on marginally trapped surfaces, including their connection to minimal surfaces, we refer the reader to the survey article \cite{AEM}.

We begin by recalling some basic definitions and properties.  By an initial data set, we mean a triple $(M,g,K)$, where $M$ is a smooth manifold, $g$ is a Riemannian metric on $M$ and $K$ is a symmetric covariant $2$-tensor on $M$.  In general relatvity, an initial data set $(M,g,K)$ corresponds to a spacelike hypersurface $M$ with induced metric $g$ and second fundamental form $K$, embedded in a spacetime (time-oriented Lorenzian manifold) $(\bar M, \bar g)$.   

Let $(M,g,K)$ be an initial data set.  For convenience, we may assume,  without loss of generality,  that this initial data set is embedded in a spacetime $(\bar M, \bar g)$ (see e.g. \cite[Section 3.2]{AM2}).   While the definition of various quantities is more natural when expressed with respect to an ambient spacetime, all the relevant quantities we introduce depend solely on the initial data set.  
With respect to the spacetime $(\bar M, \bar g)$, the tensor $K$ becomes the second fundamental form of $M$:
$K(X,Y) = \bar g(\bar\D_X u, Y)$ for all $X,Y \in T_pM$, where  $u$ is the future directed unit normal field to $M$ in 
$\bar M$.   

Let $\S$ be a closed (compact without boundary)
two-sided hypersurface in $M$. Then $\S$ admits a smooth unit normal field
$\nu$ in $M$, unique up to sign.  By convention, refer to such a choice as outward pointing. Then $\ell= u+\nu$  is a future directed outward  pointing null normal vector field along $\S$.  Associated to $\ell$  is the {\it null second fundamental form},  $\chi$  defined as, 
\beq
\chi : T_p\S \times T_p\S \to \mathbb R ,  \quad \chi(X,Y) = \bar g(\bar\D_X \ell, Y).
\eeq
In terms of the intial data,
\beq\label{chiid}
\chi = K|_{T\S} + A
\eeq
where A is the second fundamental form of $\S \subset M$ with respect to the outward unit normal $\nu$.
The {\it null expansion scalar} (or {\it null mean curvature})  $\th$ of $\S$   is obtained by tracing $\chi$
with respect to the induced metric $h$ on $\S$,
\beq
\theta = {\rm tr}_{h} \chi = h^{AB}\chi_{AB} = {\rm div}\,_{\S} \ell \,.
\eeq
Physically, $\th$ measures the divergence
of the  outgoing light rays emanating from $\S$.  In terms of the initial data $(M,g,K)$,
\beq\label{thid}
\th = {\rm tr}_{h} K + H \,,
\eeq 
where $H$ is the mean curvature of $\S$ within $M$ (given by the divergence of $\nu$ along~$\S$). 

We say that $\S$ is outer trapped (resp. weakly outer trapped) if $\th < 0$ (resp. $\th \le 0$) on $\S$.  If $\th$ vanishes identically along $\S$ then we say that $\S$ is a marginally outer trapped surface, or MOTS for short.  Note that in the so-called time-symmetric case, in which $K =0$, a MOTS is simply a minimal ($H =0$) surface in $M$, as follows from \eqref{thid}.  It is in this sense that MOTS are a spacetime generalization of minimal surfaces in Riemannian geometry.  

\subsection{Stability of MOTS.}\label{stability}
Unlike minimal surfaces, MOTS in general do not admit a variational characterization.  Nevertheless, they admit an important notion of stability which we now describe; cf., \cite{AMS2, AEM}.  
Let $\S$ be a MOTS in the initial data set $(M,g,K)$ with outward unit normal $\nu$.  Consider a normal variation of $\S$ in $M$,  i.e.,  a variation 
$t \to \S_t$ of $\S = \S_0$ with variation vector field 
$V = \frac{\d}{\d t}|_{t=0} = \phi\nu, \,\, \phi \in C^{\infty}(\S)$.
Let $\th(t)$ denote the null expansion of $\S_t$
with respect to $l_t = u + \nu_t$, where $u$ is the future
directed timelike unit normal to $M$ and $\nu_t$ is the
outer unit normal  to $\S_t$ in $M$.   A computation shows,
\beq\label{thder} \left . \frac{\d\th}{\d t} \right |_{t=0}   =
L(\f) \;, 
\eeq 
where $L : C^{\infty}(\S) \to C^{\infty}(\S)$ is
the operator~\cite{AMS2}, \beq\label{stabop2}
L(\phi)  = - \Delta \phi + 2\<X,\D\phi\>  + \left( \frac12 S_{\S}
- (\mu + J(\nu)) - \frac12 |\chi|^2+{\rm div}\, X - |X|^2
\right)\phi \,. \eeq 
In the above, $ \Delta$, $\D$ and ${\rm
div}$ are the Laplacian, gradient and divergence operator,
respectively, on $\S$, $S_{\S}$ is the scalar curvature of $\S$, $X$ is the 
vector field  on $\S$  dual to the one form $X^{\flat} = K(\nu,\cdot)|_{\S}$, 
$\<\,,\,\> =h$ is the induced metric  on $\S$, and $\mu$ and $J$ are defined in terms
of the Einstein tensor $G = {\rm Ric}_{\bar M} - \frac12 R_{\bar M} \bar g$\,: $\mu = G(u,u)$,
$J = G(u,\cdot)$.  When the Einstein equations are assume to hold, $\mu$ and $J$ represent the  energy density and linear momentum density along $M$.  As a consequence of the Gauss-Codazzi equations, the quantities $\mu$ and $J$ can be expressed solely in terms of initial
data,
\beq\label{emid}
\mu = \frac12\left(S + ({\rm tr}\,K)^2 - |K|^2\right) \quad \text{and} \quad
J = \div K- d({\rm tr}\, K)  \,,
\eeq
where $S$ is the scalar curvature on $M$.

An initial data set $(M,g, K)$ is said to satisfy the {\it dominant energy condition}, provided the inequality,
\beq\label{dec}
\mu \ge |J|
\eeq
holds along $M$.  When one assumes the Einstein equations hold, this leads to an inequality on the energy-momentum tensor that is satisfied by most classical matter fields. Note that in the time-symmetric case ($K = 0$), the dominant energy condition reduces to $S \ge 0$, and hence the importance of manifolds of nonnegative curvature in general relativity.

In the time-symmetric case, the operator $L$ reduces to the classical stability (or Jacobi) operator of minimal surface theory.   As shown in \cite{AMS2}, although $L$ is not in general self-adjoint, the eigenvalue 
$\l_1(L)$ of $L$ with the smallest real part, which is referred to as the  principal eigenvalue of $L$, is necessarily real. 
Moreover there exists an associated eigenfunction $\phi$ which is strictly positive. The MOTS 
$\S$ is then said to be stable if $\l_1(L) \ge 0$.

A basic criterion for stability is the following. We say that a  MOTS $\S$ is {\it weakly outermost} provided there are no outer trapped ($\th < 0$) surfaces outside of, and cobordant to, $\S$.  Weakly outermost MOTS are necessarily stable.  Indeed, if $\l_1(L) < 0$, Equation \eqref{thder}, with $\phi$ a positive eigenfunction ($L(\phi) = \l_1(L) \phi$)  would then imply that $\S$ could be deformed outward to an outer trapped surface.  

\subsection{Rigidity of MOTS}\label{sec:rigidity}

The proof of Theorem \ref{splitlocal1} will be based on two rigidity results for MOTS.  The following result was proved by R. Schoen and the first author in \cite{GS}.

\begin{thm}[infinitesimal rigidity]\label{pos2} 
Let  $(M,g,K)$  be an initial data set that satisfies the dominant energy condition (DEC) \eqref{dec}, $\mu \ge |J|$.
 If $\S$ is a stable MOTS in $M$ that does not admit a metric of positive scalar curvature then
\ben
\item $\S$ is Ricci flat.
\item $\chi = 0$ and $\mu + J(\nu) = 0$ along $\S$.
\een
\end{thm} 
By strengthening the stability assumption, namely by requiring the MOTS $\S$ to be  weakly outermost, as defined at the end of Section \ref{stability}, we obtain additional rigidity.
The following was proved in \cite{motsv4}.  

\begin{thm}\label{rigid2}
Let $(M,g,K)$, be an initial data set satisfying the DEC.
Suppose $\S$ is a weakly outermost MOTS in  $M$
that does not admit a metric of positive scalar curvature.  Then there exists
an outer neighborhood $U \approx [0,\delta) \times \S$ of $\S$ in $M$ such that each slice
$\S_t = \{t\} \times \S$, $t \in [0,\delta)$ is a MOTS.  
\end{thm}

\begin{rem}\label{remstable}
It follows again from the discussion at the end of Section \ref{stability} that, in the theorem above, each MOTS $\S_t$ is stable, as otherwise $\S$ would not be weakly outermost.
\end{rem}

The proofs of both rigidity results rely on the {\it MOTS stability inequality} obtained in \cite{GS}
(see Equation~2.12).  To prove Theorem \ref{splitlocal1}, we will apply these results to the initial data set  $(M,g, K = -\vep g)$.  The proof of Theorem \ref{splitlocal2}, is quite similar, where now one uses  the initial data set  $(M,g, K = \vep g)$.

\section{Proof of Theorem \ref{splitlocal1}}

\proof[Proof of Theorem \ref{splitlocal1}]  Let $(M,g)$ satisfy the assumptions of the theorem, and consider the initial data set $(M,g,K)$ where $K = -\vep g$.

We first observe that, with respect to this initial data set, the DEC \eqref{dec} holds.  Inserting 
$K = -\vep g$ into the expression for $\mu$ in \eqref{emid} leads to
\beq\label{muid}
\mu = \frac12 (S +\vep^2 n(n-1) ) =  \frac12 (S +\vep n(n-1) )\,.
\eeq
Hence, by property 1 of Theorem \ref{splitlocal1}, $\mu  \ge 0$.  Further, $K = -\vep g$ implies $J =0$, so that $\mu + |J| \ge 0$, and the DEC is satisfied. 

Next, let's consider the null expansion of $N$.  Equation \eqref{thid} implies that $N$ has null expansion,
\begin{align}\label{eqtheta}
\th &= -\vep (n-1) + H_N  \,.
\end{align}
Hence by property 2 of Theorem \ref{splitlocal1}, $\th \le 0$, i.e. $N$ is weakly outer trapped.  In fact one must have $\th \equiv 0$. Otherwise, it follows from \cite[Lemma 5.2]{AM2}, that, by a small perturbation of $N$, there would exist a strictly outer trapped ($\th < 0$) compact hypersurface $N' \subset U$ outside of, and cobordant to $N$, thereby contradicting the assumption that $N$ is weakly outermost 
in~$U$. 

Hence, $N$ is a weakly outermost MOTS in $U$.  So, by Theorem \ref{rigid2}, we can introduce coordinates $(t, x^i)$ on a neighborhood $V = [0, \delta)\times N$ of $N$ in $U$, so that $g$ in these coordinates  may be written as,
\beq\label{metric}
g = \psi^2 dt^2 + h_{ij}dx^idx^j \, ,
\eeq
where $\psi = \psi(t, x^i)$ is positive, $h_t = h_{ij}(t,x^i)dx^idx^j$ is the induced metric on $N_t = \{t\} \times N$, and $N_t$ is a MOTS, $\th(t) = 0$.

A computation similar to that leading to (\ref{thder}) (but where for the moment we do not assume 
assume $\th = \th(t)$ vanishes) leads to the following `evolution equation' for $\th = \th(t,x^i)$ (\cite{AM1, motsv3}),
\begin{align}
\frac{\d\th}{\d t}  &=  
%L_t (\phi)
 - \Delta \psi + 2\<X_t,\D\psi\> +  \left(Q_t   - \frac12 \th^2  + \th\;{\rm tr} K
+{\rm div}\, X_t - |X_t|^2 \right)\phi \,,  \label{evolve}  \\
Q_t &= \frac12 S_{N_t} - (\mu + J(\nu)) - \frac12 |\chi_t|^2\,, \label{evolve2}
\end{align}
where it is understood that, for each $t$, the above terms  live on $\S_t$,  
e.g., $ \Delta =  \Delta_t$ is the Laplacian on $N_t$, $\< , \> = h_t$, $X_t^{\flat} = K(\nu_t,\cdot)|_{N_t}$, etc.

Note from the form of $K$, $X_t = 0$.  Setting $\th =0$ and $X_t = 0$ in \eqref{evolve}, and using \eqref{evolve2},
we obtain,
\beq \label{pde}
  \Delta \psi + ( (\mu + J(\nu)) + \frac12 |\chi_t|^2 -\frac12 S_{N_t}) \psi  =0 \,.
\eeq
By Remark \ref{remstable}, each $N_t$ is a stable MOTS.  Hence, by Theorem~\ref{pos2}, 
\beq\label{vanish}
N_t \text{ is Ricci flat}, \,\, \chi_t = 0,  \text{ and } \mu + J(\nu) = 0 \,.  
\eeq
Equation \eqref{pde} then becomes,
$$
 \Delta \psi = 0 \, ,
$$
and, hence, $\psi$ is constant along each $N_t$, $\psi = \psi(t)$.  By a simple change of variable, we thus may assume $\psi = 1$, and so \eqref{metric} becomes,
\beq\label{metric2}
g = dt^2 + h_{ij}dx^idx^j \, .
\eeq

From \eqref{chiid}, $\chi_t = K|_{TN_t} +A_t = -\vep h_t + A_t$ where $A_t$ is the second fundamental form of $N_t$. 
Then, from the second equation in \eqref{vanish}, $A_t = \vep h_t$, which becomes, in the coordinate expression \eqref{metric2}, $\frac{\d h_{ij}}{\d t} = 2\vep h_{ij}$.  Integrating gives, $h_{ij}(t,x) = e^{2\vep t} h_{ij}(0,x)$.  Thus, up to isometry, we have 
$V = [0, \delta) \times N$, $g|_V = dt^2 + e^{2\varepsilon t} h$.\qed

\smallskip
\begin{rem} Theorem \ref{splitlocal1} has the following  consequence.  Let $(M,g)$  be an $n$-dimensional, $3 \le n \le 7$, asymptotically flat manifold with compact minimal boundary $N$, and with nonnegative scalar curvature, $S \ge 0$.  Suppose, further, that $N$ is an outermost minimal surface, i.e. suppose that there are no minimal surfaces in 
$M \setminus N$ homologous to $N$.  Then $N$ necessarily carries a metric of positive scalar curvature. 
For, suppose not.  To apply Theorem~\ref{splitlocal1} in the case $\vep = 0$,  it is sufficient to show that $N$ is locally weakly outermost.  If that were not the case, there would exist a compact hypersurface $N_1$ cobordant to $N$ with mean curvature $H_1 < 0$.  On the other hand sufficient far out on the asymptotically flat end there exists a compact hypersurface $N_2$ cobordant to $N_1$ with mean curvature $H_2 > 0$.  $N_1$ and $N_2$ bound a region $W$.  Basic existence results for minimal surfaces (or for MOTS \cite{AEM}), guarantee the existence of a minimal surface in $W$ homologous to $N$, contrary to assumption.  Hence $N$ is weakly outermost.  Theorem \ref{splitlocal1} then implies that $(M,g)$ locally splits near $N$, contrary to $N$ being an outermost minimal surface. The same consequence holds for an $n$ dimensional, $3 \le n \le 7$, asymptotically hyperbolic manifold $(M,g)$ with compact boundary $N$ of constant mean curvature $n-1$, and with scalar curvature $S\ge -n(n-1)$ in the following sense: Consider the initial data set $(M,g,-g)$, with $(M,g)$ as just described, and suppose $N$ is an outermost MOTS. Then $N$ necessarily carries a metric of positive scalar curvature.
\end{rem}

\smallskip
Theorem \ref{splitlocal1} globalizes in a straight-forward way, as follows.

\begin{thm}\label{split}
Let $(M,g)$ be a complete, noncompact $n$-dimensional ($n \ge 3$) Riemannian manifold with compact boundary $N$.  Assume:
\ben
\item $M$ has scalar curvature $S \ge  -\varepsilon n(n-1)$, where $\varepsilon = 0$ or $1$.
\item $N$ has mean curvature $H_N \le  \varepsilon (n-1)$.
\item $N$ does not carry a metric of positive scalar curvature and is weakly outermost.
\een
Then $(M,g)$ is isometric to $[0,\infty) \times N$, with (warped) product metric $dt^2 + e^{2\varepsilon t} h$, 
where $(N,h)$ is Ricci flat. 
\end{thm}

\proof[Proof of Theorem \ref{split}]  By Theorem \ref{splitlocal1}, there exists a neighborhood $V$ of $N$ such that $(V, g|_V)$ is isometric to $([0, \delta) \times N, dt^2 + e^{2\varepsilon t} h)$.  By the completeness assumption, it is clear that this warped product structure extends to $t =\delta$. From the fact that $N$ is weakly outermost, it follows that $N_{\delta} = \{\delta\} \times N$ is  weakly outermost.  Theorem \ref{splitlocal1} then implies that the warped product structure extends beyond $t = \delta$.  By a continuation argument, it follows that the warped product structure exists for all $t \in [0,\infty)$.\qed

\proof[Proof of Theorem \ref{splitlocal2}]  The proof of Theorem \ref{splitlocal2} is very similar to the proof of  \ref{splitlocal1}, except that now one works with the initial data set 
$(M,g, K = \vep g)$. We leave the details to the interested reader.

\smallskip
Similar to Theorem \ref{splitlocal1}, Theorem \ref{splitlocal2} implies  the following global result.

\begin{thm}\label{split2}
Let $(M,g)$ be a complete, noncompact $n$-dimensional ($n \ge 3$) Riemannian manifold with compact boundary $N$.  Assume:
\ben
\item $M$ has scalar curvature $S \ge  -\varepsilon n(n-1)$, where $\varepsilon = 0$ or $1$.
\item $N$ has mean curvature $H_N \le - \varepsilon (n-1)$.
\item $N$ does not carry a metric of positive scalar curvature and is weakly outermost.
\een
Then $(M,g)$ is isometric to $[0,\infty) \times N$, with (warped) product metric $dt^2 + e^{-2\varepsilon t} h$, 
where $(N,h)$ is Ricci flat. 
\end{thm}

\section{Warped product splitting and Obata's equation}

The main aim of this section is to prove Theorem \ref{splithess} stated in the introduction.
  Obata's equation in the form $\nabla^2 f=fg$ has been studied previously in the literature; see e.g. \cite{Kanai, Tashiro}. In addition to Theorems  \ref{splitlocal1} and \ref{splitlocal2}, the proof of 
Theorem~\ref{splithess} will make use of the following result, which extends to manifolds with boundary certain results in \cite{Kanai}.

\begin{prop}\label{proposition:rigidity by hessian eqn}
Let $(M^n,g)$ be a complete connected Riemannian manifold with compact connected boundary $ N$ $(n\ge 3)$. Suppose there exists a nonzero function $f$ that satisfies 
\beq\label{eq:hess}
\nabla^2 f=fg,
\eeq
and $ N$ is a regular hypersurface $f^{-1}(a)$ for $a\in\mathbb{R}$. Then the following hold:
\begin{enumerate}
	\item If $M$ is compact, then $(M,g)$ is isometric to a hyperbolic cap $[0,R]\times \mathbb{S}^{n-1}$ equipped with the metric
	\[
		dt^2+(\sinh t)^2 g_{\mathbb{S}^{n-1}}
	\]
	where $g_{\mathbb{S}^{n-1}}$ is the standard unit sphere metric and $R=d_g(p, N)$ for $p\in M\setminus N$ which is a critical point of $f$.
	\item If $M$ is noncompact, then $(M,g)$ is isometric to a manifold $[0,\infty)\times N$ with (warped) product metric of the form
		\[dt^2+\xi(t)^2 g|_ N,\]
	where  $\xi:[0,\infty)\rightarrow\mathbb{R}$ is the solution to the following ODE
	\begin{equation}\label{equation:warping factor}
		\left\{
		\begin{aligned}
			&\xi''-\xi=0 \text{ on }[0,\infty),\\
			&\xi(0)=1 \text{ and }\xi'(0)=\frac{a}{|\nabla f|_N} \,.		
			\end{aligned}\right.
	\end{equation}
(We note, as follows from \eqref{eq:hess}, that $|\nabla f|_N$ is constant.)
\end{enumerate}
\end{prop}

\proof[Proof of Proposition \ref{proposition:rigidity by hessian eqn}] First we claim that $f$ has a critical point on the interior of $M$ if and only if $M$ is compact (with boundary).

	Suppose that $f$ has a critical point $p$ in $M$. Consider a unit speed geodesic $\gamma:[0,\infty)\to M$ emanating from $p$. It follows that
	\[
		\frac{d^2}{dr^2}f(\gamma(r))-f(\gamma(r))=0
	\]
thus $f(\gamma(r))=c(e^r+e^{-r})$ and $\frac{d}{dr}f(\gamma(r))=c(e^r-e^{-r})$, where $c\ne 0$
(as otherwise $f$ would vanish identically). Observe that $f$ depends only on the geodesic distance from the point $p$, which implies that $\gamma'$ is parallel to $\nabla f$. Moreover, there cannot be any other critical point of $f$. Let $R=d_g(p, N)$. Then it follows that $ N=\text{exp}_p(S_R)$, and from this that $\text{exp}_p:\overline{B_R}\to M$ is bijective. By continuity of the exponential map, this implies that $M$ must be compact.

	Suppose, conversely, $M$ is compact. For contradiction, suppose also that $f$ has no critical points. Without loss of generality, we may assume that $\nabla f$ points inward on $N$. Let 
$\nu=\nabla f/|\nabla f|$, 
and consider the integral curve $\gamma$ of $\nu$ emanating from a point $p\in N$, i.e., $\gamma(0)=p$. It is straightforward that $\gamma$ is a geodesic parametrized by arc length, and we also have
	\[
		f\circ\gamma(t)=c_1e^t+c_2e^{-t}
	\]
	as we observed before. Since $f$ has no critical point, $\gamma$ can be extended to $[0,\infty)$, which implies that $\gamma$ is an injective infinite length  geodesic. This contradicts the condition that $M$ is compact, hence $f$ must have a critical point on the interior of $M$.

	We now show the first case of the proposition: assume that $M$ is compact. From the previous argument, there is a critical point $p$ such that $\text{exp}_p:\overline{B_R}\to M$ is bijective where $R=d_g(p, N)$. Now we show that it is a diffeomorphism. Let $J$ be a Jacobi field along $\gamma$ such that $J(0)=0$ and $|J'(0)|=1$ and $g(J',\gamma')=0$. Then we have for $r>0$,
		\begin{align*}
		\left.\frac{1}{2}\frac{d}{dt}\right|_{t=r} g(J,J)&=g(J,\nabla_{\gamma'}J)|_{t=r}\\
		&=A(J,J)|_{t=r}=\left.\frac{f}{|\nabla f|}g(J,J)\right|_{t=r} \, ,
	\end{align*}
where $A = \nabla^2 f/|\D f|$ is the second fundamental form of the geodesic spheres.
Thus for $r>r_0>0$,
	\begin{align*}
		|J|^2(r)=\left(\frac{ e^r- e^{-r}}{e^{r_0}- e^{-r_0}}\right)^2|J|^2(r_0)\ne 0
	\end{align*}
where $r_0$ is sufficiently small that $|J(r_0)|\ne 0$. This implies that there is no conjugate point from $p$ thus $\text{exp}_p:\overline{B_R}\to M$ is a diffeomorphism. Furthermore, by using geodesic polar coordinates, we can write the metric $g$ on $M$ diffeomorphic to $[0,R]\times\mathbb{S}^{n-1}$ as
	\[
		g=dt^2+(\sinh t)^2 g_{\mathbb{S}^{n-1}}.
	\]
	
	We turn to the second case: assume that $M$ is noncompact. Let $h=g|_ N$. We will construct an isometry between $(M,g)$ and the manifold 
	\[([0,\infty)\times  N, dt^2+\xi^2 h)\]
	where $\xi$ is given in (\ref{equation:warping factor}). Without loss of generality, we may assume that $\nabla f$ points inward on $ N$. 

Let $\varphi$ be the flow generated by $\nu=\nabla f/|\nabla f|$, and define the map $\psi:[0,\infty)\times N\to M$ by 
	\[\psi(t,\bar{p})=\varphi_{t}(\bar{p})\] for $\bar{p}\in  N$ and $t\in[0,\infty)$.
	Since $f$ has no critical points, it is clear that $\psi$ is a diffeomorphism. 
	
As we observed before, we have the general solution
	\[f\circ\psi(\bar{p},t)=c_1e^{t}+c_2e^{-t},\quad\bar{p}\in  N,t\in[0,\infty) \, ,\]  
where the  constants $c_1$ and $c_2$ are determined by the conditions on $f$ at $N$; specifically, $c_1 +c_2 = a$ and $c_1 - c_2 = |\D f|_N$. In terms of these constants, the solution to the ODE \eqref{equation:warping factor} is given by, $\xi(t)=\frac{c_1e^{t}-c_2e^{-t}}{c_1-c_2}$.

Now we prove that $\psi$ is the desired isometry from $([0,\infty)\times N,dt^2+\xi(t)^2 h)$ onto $(M,g)$. Using $\psi$ as a coordinate chart, we can write the metric
	\[
		g=dt^2+g_{ij}(t,\bar{p})\,dx^idx^j
	\] where $\{x^i\}_{i=1}^{n-1}$ are local coordinates near $\bar{p}$ on $ N$ and $g_{ij}(t,\bar{p})=g_{(t,\bar{p})}(\partial_i,\partial_j)$ for $t\in [0,\infty)$. Then we have
	\begin{align*}
		\left.\frac{1}{2}\frac{d}{dt}\right|_{t=\tau}g_{ij}(t,\bar{p})&=g_{(\tau,\bar{p})}(\partial_i,\nabla_\nu \partial_j)=A_{(\tau,\bar{p})}(\partial_i,\partial_j)\\
		&=\frac{f(\tau)}{|\nabla f|(\tau)}g_{ij}(\tau,\bar{p})=\frac{c_1e^{\tau}+c_2e^{-\tau}}{c_1e^{\tau}-c_2e^{-\tau}}g_{ij}(\tau,\bar{p})
	\end{align*}
	where $A_{(\tau,\bar{p})}$ is the second fundamental form of the hypersurface $f^{-1}(f\circ\gamma(\tau))$. Thus we obtain 
	\[
		g_{ij}(\tau,\bar{p})=\xi(\tau)^2 g_{ij}(0,\bar{p})=\xi(\tau)^2 h_{ij}(\bar{p})  \,,
	\] and by varing $(\tau,\bar{p})\in M$ it proves that $\psi$ is the desired isometry.\qed

\smallskip
\proof[Proof of Theorem \ref{splithess}] We will only prove $\delta=1$ since the proof of the other case is almost identical.

 By Theorem \ref{splitlocal1} (with $\vep = 1$), we have the local splitting near $ N$, that is, there exists a neighborhood $U$ of $ N$ such that $U$ is isometric to $[0,b)\times N$ for some $b>0$ with the metric $dt^2+e^{2t} h$.
 To use Proposition \ref{proposition:rigidity by hessian eqn}, we show that $ N$ is the level set $f^{-1}(a)$ for some $a\in\mathbb{R}$.

Let $\{x^i\}_{i=1}^{n-1}$ be local coordinates on $N$. This gives rise to local coordinates 
$\{t = x_0, x_1,...,x_{n-1}\}$ on  $[0,b)\times N$ in the obvious manner.   
  Then, by direct computation, we have
	\begin{align}
		0=\nabla_{\partial_t}\nabla_{\partial_i}f&=\partial_t\partial_i f-\sum_{k=0}^{n-1}\Gamma_{ti}^k\partial_k f\nonumber\\
		&=\partial_t\partial_i f-\partial_i f,\\
		f=\nabla_{\partial_t}\nabla_{\partial_t}f&=\partial^2_t f,\\
		e^{2t}f\, h_{ij}=\nabla_{\partial_i}\nabla_{\partial_j}f&=\partial_i\partial_j f+e^{2t}h_{ij}\partial_t f-\sum_{l=1}^{n-1}\bar{\Gamma}_{ij}^{l}\partial_l f.\label{eqn:level set eqn0}
	\end{align}
	where $\bar{\Gamma}$ is the Christoffel symbol with respect to $h$.
	Denote $f=f(p,t)$ on $U$ for $p\in  N$ and $t\in [0,b)$. Then from the above computations we have
	\begin{align}
		\partial_t^2f-f=0 &\Rightarrow f(p,t)=c_1(p)e^t+c_2(p)e^{-t},\label{eqn:level set eqn1}\\
		\partial_t(\partial_i f)-\partial_i f=0 &\Rightarrow \partial_i f(p,t)=c_3(p)e^t,\label{eqn:level set eqn2}
	\end{align}
	It follows from \eqref{eqn:level set eqn1} and \eqref{eqn:level set eqn2} that $c_2(p)$ is constant on $ N$, hence we can write
	\[
		f(p,t)=c_1(p)e^t+c_2e^{-t},\text{ and } \partial_t f-f=-2c_2e^{-t}.
	\]
	Now we shall show that $c_1(p)$ is constant on $ N$. By \eqref{eqn:level set eqn0}, we have
	\begin{align}
		\partial_i\partial_j f+e^{2t}h_{ij}&(\partial_t f-f)-\sum_{l=1}^{n-1}\bar{\Gamma}_{ij}^l\partial_l f=0\nonumber\\
		&\Rightarrow e^t\left(\partial_i\partial_j (c_1(p))-2c_2 h_{ij}-\sum_{l=1}^{n-1}\bar{\Gamma}_{ij}^{l}\partial_l c_1(p)\right)=0\nonumber\\
		&\Rightarrow \nabla^ N_{\partial_i}\nabla^ N_{\partial_j}c_1=2c_2 h_{ij} \Rightarrow  \Delta_ N c_1=2(n-1)c_2  \,. \label{eqn:level set eqn3}
	\end{align}
	Since $ N$ is compact without boundary, we have
	\[
		0=\int_ N  \Delta_ N c_1=2(n-1)c_2| N|
	\]
	where $| N|$ is the area of $ N$. This implies that $c_2=0$, and hence $c_1(p)$ is harmonic on $ N$. Therefore $c_1(p)$ is constant on $ N$ so $ N=f^{-1}(c_1)$.

	By Proposition \ref{proposition:rigidity by hessian eqn}, $(M,g)$ is isometric to $[0,\infty)\times N$ with the metric $dt^2+\xi(t)^2 h$. In particular, one can see that the warping factor is $\xi(t)=e^t$. \qed

\bigskip
\noindent
\textsc{Acknowledgements.}   GJG was partially supported by  NSF grant DMS-1710808. HCJ was partially supported by NSF Grant DMS-1452477.  The authors would like to thank Lan-Hsuan Huang for her interest in this work and for many helpful comments. HCJ is especially grateful to Professor Huang for her constant encouragement and guidance. We would also like to thank an anonymous referee for some helpful comments.

\bibliographystyle{amsplain}
\bibliography{splitting}

\end{document}